\documentclass{amsart}
\usepackage{amssymb}

\usepackage{graphicx}
\usepackage{amscd}

\newtheorem{theorem}{Theorem}
\theoremstyle{plain}

\newtheorem{corollary}{Corollary}

\newtheorem{lemma}{Lemma}

\newtheorem{proposition}{Proposition}
\newtheorem{remark}{Remark}

\numberwithin{equation}{section}
\numberwithin{theorem}{section}
\numberwithin{lemma}{section}
\numberwithin{proposition}{section}
\numberwithin{corollary}{section}
\numberwithin{remark}{section}
\numberwithin{example}{section}
\numberwithin{definition}{section}

\input{tcilatex}

\begin{document}
\title[ Local Zeta Functions and Linear Feedback Shift Registers]{Computing \ Igusa's Local Zeta Functions of Univariate Polynomials, and
Linear Feedback Shift Registers}
\author{W. A. Zuniga-Galindo}

\begin{abstract}
We give a polynomial time algorithm for computing the Igusa local zeta
function $Z(s,f)$ attached to a polynomial \ $f(x)\in \mathbb{Z}[x]$, in one
variable, with splitting field $\mathbb{Q}$, and a prime number $p$. We also
propose a new class of Linear Feedback Shift Registers based on the
computation of Igusa's local zeta function.
\end{abstract}

\address{Department of Mathematics and Computer Science, Barry University, 11300 N.E.
Second Avenue\\
Miami Shores, Florida 33161, USA}
\email{wzuniga@mail.barry.edu}
\thanks{Supported by COLCIENCIAS-Grant \# 089-2000. }
\maketitle

\section{Introduction}

Let $f(x)\in \mathbb{Z}[x]$, $x=(x_{1},\mathbb{\cdots },x_{n})$ be a
non-constant polynomial, and $\ p$ a fixed prime number. We put $%
N_{m}(f,p)=N_{m}(f)$ for the number of solutions of the congruence $%
f(x)\equiv 0$ mod $p^{m}$ in $(\mathbb{Z}/p^{m}\mathbb{Z})^{n}$, $m\geqq 1$,
\ and \ $H(t,f)$ for the Poincar\'{e} series 
\begin{equation*}
H(t,f)=\sum\limits_{m=0}^{\infty }N_{m}(f)(p^{-n}t)^{m}\text{,}
\end{equation*}
with $t\in \mathbb{C}$, $\mid t\mid <1$, and \ $N_{0}(f)=1$. This paper is
dedicated to the computation of the sequence $\left\{ N_{m}(f)\right\}
_{m\geqq 0}$ when $f$ is an univariate polynomial with splitting field $%
\mathbb{Q}$.

Igusa showed that the \ Poincar\'{e} series $H(t,f)$ admits a meromorphic
continuation to the complex plane as a rational function of $t$ \cite{I1}, 
\cite{I2}. In this paper we make a first step towards the solution of the
following problem: given a polynomial \ $f(x)$\ as above, how difficult is
to compute the meromorphic continuation of the \ Poincar\'{e} series $H(t,f)$%
?

The computation of the \ Poincar\'{e} series $H(t,f)$\ is equivalent to the
computation of Igusa's local zeta \ function $Z(s,f)$, attached to $f$ and $%
p $, defined as follows. We denote by $\mathbb{Q}_{p}$ the field of $p-$adic
numbers, and by $\mathbb{Z}_{p}$ the ring of $p-$adic integers. For $x\in 
\mathbb{Q}_{p}$, $v_{p}(x)$ denotes the $p-$adic order of $x$, and $%
|x|_{p}=p^{-v_{p}(x)}$ its absolute value. The Igusa local zeta function
associated to $f$ and $p$ is defined as follows: 
\begin{equation*}
Z(s,f)=\int_{\mathbb{Z}_{p}^{n}}|f(x)|_{p}^{s}\mid dx\mid
,\,\,\,\,\,\,\,\,s\in \mathbb{C},\,\,
\end{equation*}
where $\func{Re}(s)>0$, and $\mid dx\mid $ denotes the Haar measure on $%
\mathbb{Q}_{p}^{n}$ so normalized that $\mathbb{Z}_{p}^{n}$ has measure $1$.
The following relation between $Z(s,f)$\ and $\ H(t,f)$ holds (see \cite{I1}%
, theorem 8.2.2): 
\begin{equation*}
H(t,f)=\frac{1-tZ(s,f)}{1-t},\text{ }t=p^{-s}.
\end{equation*}
Thus, the rationality of $Z(s,f)$\ implies the rationality of the
Poincar\'{e} series \ $H(t,f)$, and \ the computation of $H(t,f)$ is \
equivalent to the computation of$\ Z(s,f)$. Igusa \cite[theorem 8.2.1]{I1}\
showed that the local zeta function $Z(s,f)$ \ admits a meromorphic
continuation to the complex plane as a rational function \ of $p^{-s}$.

The first result of this paper is a polynomial time algorithm for computing
the local zeta function $Z(s,f)$ attached to a polynomial \ $f(x)\in \mathbb{%
Z}[x]$, in one variable, with splitting field $\mathbb{Q}$, and a prime
number $p$. We also give an explicit estimate for its complexity (see
algorithm Compute\_$Z(s,f)$ in section 2, and theorem \ref{theoA}).

Many authors have found explicit formulas for $Z(s,f)$, or $H(f,t)$,\ for
several classes of polynomials, among them \cite{D2}, \cite{D-H}, \cite{G1}, 
\cite{G2}, [\cite{I3} and the references therein], \cite{S-Z-G}, \cite{Z-G1}%
, \cite{Z-G2}. In all these works the computation of $Z(s,f)$, or $H(f,t)$,
is reduced to the computation of \ other problems, as the computation of the
number of solutions of polynomial equations with coefficients in a finite
field. Currently, there is no polynomial time algorithm solving this problem
\ \cite{Wa}, \cite{V-K-S}. Moreover, none of the above mentioned works
include complexity estimates for the computation of Igusa's local zeta
functions.

Of particular importance is Denef's explicit formula for \ $Z(s,f)$, when $f$
satisfies some generic conditions \ \cite{D2}. This formula involves the
numerical data associated to a resolution of singularities of the divisor $%
f=0$, and the number of rational points of certain non-singular varieties
over finite fields. Thus the computation of $Z(s,f),$ for a generic
polynomial $f$, is reduced to the computation of \ the numerical data
associated to a resolution of singularities of the divisor $f=0$, and the
number of solutions of non-singular polynomials over finite fields.
Currently, it is unknown if these problems can be solved in polynomial time
on a Turing machine. However, during the last few years important
achievements have been obtained in the computation of resolution of
singularities of polynomials \cite{B-M}, \cite{B-S}, \cite{B-S1}, \cite{Vi}.

The computation of the Igusa local zeta function for an arbitrary polynomial
seems to be an intractable problem on a Turing machine. For example, for $%
p=2 $, the computation of the number of solutions of a polynomial equation
with coefficients in $\mathbb{Z}/2\mathbb{Z}$ is an $\mathbf{NP}-$complete
problem on a Turing Machine \cite[page 251, problem AN9]{G-J}. Then \ in the
case of $2-$adic numbers, the computation of \ the Igusa local zeta function
is an $\mathbf{NP}-$complete problem.

Recently, Anshel and Goldfeld have shown the existence of a strong
connection between the computation of zeta functions and cryptography \cite
{A-G}. Indeed, they proposed a new class of candidates for one-way functions
based on global zeta functions. \ A one-way function is a function $\ F$
such that for each $x$ in the domain of $F$, it is easy to compute $F(x)$;
but for essentially all $y$ in the range of $F$, it is an intractable
problem to find an $x$ \ such that $y=F(x)$. These functions play a central
role, from a practical and theoretical point of view, in modern
cryptography. Currently, there is no guarantee that one-way functions exist
even if $\mathbf{P}\neq \mathbf{NP}$. Most of the present candidates for
one-way functions are constructed on the intractability of problems like
integer factorization and discrete logarithms \cite{G-L-N}. Recently, P.
Shor has introduced a new approach to attack these problems \cite{S}.
Indeed, Shor have shown that on a quantum computer the integer factorization
and discrete logarithm \ problems can be computed in polynomial time.

We set 
\begin{equation*}
\mathcal{H}=\{H(t,f)\mid f(x)\in \mathbb{Z}[x]\text{, in \ one variable,
with splitting field }\mathbb{Q}\},
\end{equation*}
and $N^{\infty }$ $\left( \mathbb{Z}\right) $ for the set \ of finite
sequences of integers. For each positive integer $u$ and a prime number $p$,
we define 
\begin{equation*}
\begin{array}{cccc}
F_{u,p}: & \mathcal{H} & {\ \rightarrow } & \mathbb{N}^{\infty }\left( 
\mathbb{Z}\right) \\ 
& H(t,f) & {\ \rightarrow } & \{N_{0}(f,p),N_{1}(f,p),\cdots ,N_{u}(f,p)\}%
\text{.}
\end{array}
\end{equation*}

Our second result asserts that \ $F_{u,p}(H(t,f))$ \ can be computed in
polynomial time, for every $H(t,f)$\ in $\mathcal{H}$ (see theorem \ref
{theoB}). It seems interesting \ to study the complexity on a Turing machine
of the following \ problem: given a list of positive integers $%
\{a_{0},a_{1},\cdots ,a_{u}\}$, how difficult is it to determine whether or
not there exists a Poincar\'{e} series $H(t,f)=\sum_{m=0}^{\infty
}N_{m}(f)(p^{-1}t)^{m}$, such that $a_{i}=N_{i}(f)$, $i=1$,$\cdots $, $u$?

Currently, the author does not have any result about the complexity of the
above problem, however the mappings $F_{u,p}$ can be considered as new \
class of stream ciphers (see section 8).

\section{The Algorithm \ Compute\_$Z(s,f)$}

In this section we present a polynomial time algorithm, Compute\_$Z(s,f)$, \
that solves the following problem: given \ a polynomial $f(x)\in \mathbb{Z}%
[x]$, in one variable, whose splitting field is $\mathbb{Q}$, \ find an
explicit expression for the meromorphic continuation of $Z(s,f)$. The
algorithm is as follows.

\textbf{Algorithm\ Compute\_}$Z(s,f)$

Input\ : A \ polynomial $f(x)\in \mathbb{Z}[x],$\ in one variable, whose
splitting field is $\mathbb{Q}$.

Output : A rational function of $p^{-s}$\ that is the meromorphic \qquad
continuation of $Z(s,f)$.

\begin{enumerate}
\item  Factorize $f(x)$ in $\mathbb{Q}[x]$: $f(x)=\alpha
_{0}\prod\limits_{i=1}^{r}(x-\alpha _{i})^{e_{i}}\in \mathbb{Q}[x].\qquad $

\item  Compute 
\begin{equation*}
l_{f}=\left\{ 
\begin{array}{cc}
1+\max \{v_{p}(\alpha _{i}-\alpha _{j})\mid i\neq j\text{, }1\leq i\text{, }%
j\leq r\}\text{,} & \text{ if }r\geqq 2\text{;} \\ 
1\text{, } & \text{\ if \ \ }r=1\text{.}
\end{array}
\right. 
\end{equation*}

\item  Compute the $p-$adic expansions of the numbers $\alpha _{i}$, $%
i=1,2,\cdots ,r$\ modulo $p^{l_{f}+1}$.

\item  Compute the tree $T(f,l_{f})$\ associated to $f(x)$\ and $p$ (for the
definition of $T(f,l_{f})$ see (\ref{tree})).

\item  Compute the generating function $G(s,T(f,l_{f}),p)$\ attached to $%
T(f,l_{f})$ (for the definition of $G(s,T(f,l_{f}),p)$\ see (\ref{Gfunction}%
)).

\item  Return $Z(s,f)=G(s,T(f,l_{f}),p)$.

\item  End
\end{enumerate}

In section 6, we shall give a proof of the correctness \ and a \ complexity
estimate for the algorithm Compute\_$Z(s,f)$. The first step in our
algorithm is accomplished by means \ of the factoring algorithm by A.K.
Lenstra, H. Lenstra and L. Lov\'{a}sz \cite{L-L-L}. If $d_{f}$ denotes the
degree of $f(x)=\sum\nolimits_{i}a_{i}x^{i}$, and 
\begin{equation*}
\shortparallel f\shortparallel =\sqrt{\sum_{i}a_{i}^{2}},
\end{equation*}
then\ \ the \ mentioned factoring algorithm needs \ $O\left(
d_{f}^{6}+d_{f}^{9}(\log \shortparallel f\shortparallel )\right) $
arithmetic operations, and the integers on which these operations are
performed each have a binary length \ \ $O\left( d_{f}^{3}+d_{f}^{2}(\log
\shortparallel f\shortparallel )\right) $ \cite[theorem 3.6]{L-L-L}.

The steps $2$, $3$, $4$, $5$ reduce \ in polynomial time the computation of $%
Z(s,f)$ to the computation of a factorization of $f(x)$ over $\mathbb{Q}$.
This reduction is accomplished by constructing a weighted tree from the $p-$%
adic expansion of the roots of $\ f(x)$ modulo a \ certain power of $p$ (see
section 4), and then \ \ associating a generating function to this tree (see
section 5). Finally, \ we shall prove that the generating function \
constructed in this way coincides with the local zeta function of $f(x)$
(see section 5).

\section{$p$-adic Stationary Phase Formula}

Our main tool in the effective computing of Igusa's local zeta function \ of
a polynomial in one variable will be the $p-$adic stationary phase formula,
abbreviated SPF \cite{I3}. This formula is a recursive procedure for
computing local zeta functions. By using this procedure it is possible to
compute the local zeta functions for many classes of polynomials [\cite{I3}
and the references therein], \cite{S-Z-G}, \cite{Z-G1}, \cite{Z-G2}, \cite
{ZG3}.

Given a polynomial $f(x)\in $ $\mathbb{Z}_{p}[x]\setminus $ $p\mathbb{Z}%
_{p}[x]$, we denote by $\overline{f(x)}$ its reduction modulo $p\mathbb{Z}%
_{p}$, i.e., the polynomial obtained by reducing the coefficients of $f(x)$
modulo $p\mathbb{Z}_{p}$. \ We define for each $x_{0}\in $ $\mathbb{Z}_{p}$, 
\begin{equation*}
f_{x_{0}}(x)=p^{-e_{x_{0}}}f(x_{0}+px)\text{,}
\end{equation*}
where $e_{x_{0}}$\ is the minimum \ order of $p$\ \ \ in the coefficients of 
$f(x_{0}+px)$. Thus $f_{x_{0}}(x)\in \mathbb{Z}_{p}[x]\setminus $ $p\mathbb{Z%
}_{p}[x]$. We shall call the polynomial $f_{x_{0}}(x)$\ the \textit{%
dilatation} of $f(x)$\ at $x_{0}$. We also define 
\begin{equation*}
\nu (\overline{f})=\text{Card}\{\overline{z}\in \mathbb{F}_{p}\mid \overline{%
f}(\overline{z})\neq 0\},
\end{equation*}
\begin{equation*}
\delta (\overline{f})=\text{Card}\{\overline{z}\in \mathbb{F}_{p}\mid 
\overline{z}\text{ is a simple root of }\overline{f}(\overline{z})=0\}.
\end{equation*}

We shall use $\{0,1,\cdots ,p-1\}\subseteq \mathbb{Z}_{p}$ as a set of
representatives of the elements of $\mathbb{F}_{p}=\mathbb{Z}/$ $p\mathbb{Z}%
=\{\overline{0},\overline{1},\cdots ,\overline{p-1}\}$. Let $S=S(f)$ denote
the subset of $\{0,1,\cdots ,p-1\}\subseteq \mathbb{Z}_{p}$ which is mapped
bijectively \ by the canonical homomorphism $\mathbb{Z}_{p}\rightarrow 
\mathbb{Z}_{p}/p\mathbb{Z}_{p}$ to the set of roots of \ $\overline{f}(%
\overline{z})=0$ with multiplicity greater than or equal to two.

With all the above notation we are able to state the $p-$adic stationary
phase formula for polynomials in one variable.

\begin{proposition}[{\protect\cite[theorem 10.2.1]{I1}}]
Let $f(x)\in $ $\mathbb{Z}_{p}[x]\setminus $ $p\mathbb{Z}_{p}[x]$ be a
non-cons\-tant polynomial. Then 
\begin{equation*}
Z(s,f)=p^{-1}\nu (\overline{f})+\delta (\overline{f})\frac{(1-p^{-1})p^{-1-s}%
}{(1-p^{-1-s})}+\sum\limits_{\xi \in S}p^{-1-e_{\xi }s}\int\limits_{\mathbb{Z%
}_{p}}\mid f_{\xi }(x)\mid _{p}^{s}dx.
\end{equation*}
\end{proposition}

The following example illustrates the use \ of the $p-$adic stationary phase
formula, and \ also the \ basic aspects of our algorithm for computing $%
Z(s,f)$.

\subsection{Example\label{ex2}}

Let \ $f(x)=(x-\alpha _{1})(x-\alpha _{2})^{3}(x-\alpha _{3})(x-\alpha
_{4})^{2}(x-\alpha _{5})$ be a polynomial such that $\alpha _{1}$, $\alpha
_{2}$ , $\alpha _{3}$, $\alpha _{4}$ , $\alpha _{5}$ are integers having the
following $\ p-$adic expansions: 
\begin{equation*}
\alpha _{1}=a+dp+kp^{2},
\end{equation*}
\begin{equation*}
\alpha _{2}=a+dp+lp^{2},
\end{equation*}
\begin{equation*}
\alpha _{3}=b+gp+mp^{2},
\end{equation*}
\begin{equation*}
\alpha _{4}=c+hp+np^{2},
\end{equation*}
\begin{equation*}
\alpha _{5}=c+hp+rp^{2},
\end{equation*}
where \ the $p$-adic digits $a$, $b$, $c$, $d$, $g$, $h$, $l$, $m$, $n$, $r$
belong to $\{0,1,\mathbb{\cdots },p-1\}$. We assume the $p$-adic digits to
be different by pairs. The local zeta function $Z(s,f)$ will be computed by
using SPF iteratively.

By applying SPF with $\overline{f(x)}=(x-\overline{a})^{4}(x-\overline{b})(x-%
\overline{c})^{3}$, $\nu (\overline{f})=p-3$, $\delta (\overline{f}%
)=1,S=\{a,c\}$, $f_{a}(x)=p^{-4}f(a+px)$, and $f_{c}(x)=p^{-3}f(c+px)$, we
obtain that

\begin{eqnarray}
Z(s,f) &=&p^{-1}(p-3)+\frac{(1-p^{-1})p^{-1-s}}{1-p^{-1-s}}+p^{-1-4s}\int_{%
\mathbb{Z}_{p}}|f_{a}(x)|_{p}^{s}\mid dx\mid  \notag \\
&&+p^{-1-3s}\int_{\mathbb{Z}_{p}}|f_{c}(x)|_{p}^{s}\mid dx\mid \text{.}
\label{step8}
\end{eqnarray}
We apply SPF to the integrals involving $f_{a}(x)$ and $f_{c}(x)$\ in (\ref
{step8}). First, we consider the integral corresponding to \ $f_{a}(x)$.
Since $\overline{f_{a}(x)}=(x-\overline{d})^{4}(\overline{a}-\overline{b})(%
\overline{a}-\overline{c})^{3}$, $S=\{d\}$, $f_{a,d}(x)=p^{-4}f_{a}(d+px)$, $%
\nu (\overline{f_{a}})=p-1$, and $\delta (\overline{f_{a}})=0$, it follows
from (\ref{step8}) using SPF that

\begin{eqnarray}
Z(s,f) &=&p^{-1}(p-3)+\frac{(1-p^{-1})p^{-1-s}}{1-p^{-1-s}}%
+p^{-1}(p-1)p^{-1-4s}  \notag \\
&&+p^{-2-8s}\int_{\mathbb{Z}_{p}}|f_{a,d}(x)|_{p}^{s}\mid dx\mid
+p^{-1-3s}\int_{\mathbb{Z}_{p}}|f_{c}(x)|_{p}^{s}\mid dx\mid \text{.}
\label{step9}
\end{eqnarray}
Now, we apply SPF to the integral involving $f_{c}(x)$ in (\ref{step9}).
Since $\overline{f_{c}(x)}$ $=(\overline{c}-\overline{a})^{4}(\overline{c}-%
\overline{b})(x-\overline{h})^{3}$, $S=\{h\},f_{c,h}(x)=p^{-3}f_{c}(h+px)$, $%
\nu (\overline{f_{c}})=p-1$, and $\delta (\overline{f_{c}})=0$, it follows
from (\ref{step9}) using SPF that 
\begin{eqnarray}
Z(s,f) &=&p^{-1}(p-3)+\frac{(1-p^{-1})p^{-1-s}}{1-p^{-1-s}}%
+p^{-1}(p-1)p^{-1-4s}  \notag \\
&&\,+p^{-2-8s}\int_{\mathbb{Z}_{p}}|f_{a,d}(x)|_{p}^{s}\mid dx\mid
+p^{-1}(p-1)p^{-1-3s}  \notag \\
&&\,+p^{-2-6s}\int_{\mathbb{Z}_{p}}|f_{c,h}(x)|_{p}^{s}\mid dx\mid \text{.}
\label{step10}
\end{eqnarray}

By applying SPF to the integral involving $f_{a,d}(x)$ in (\ref{step10}),
with $\overline{f_{a,d}(x)}=(x-\overline{k})(x-\overline{l})^{3}(\overline{d}%
-\overline{b})(\overline{d}-\overline{c})^{3}$, $S=\{k,l\}$, $%
f_{a,d,k}(x)=p^{-1}f_{a,d}(k+px)$, $|f_{a,d,k}(x)|_{p}^{s}=|x|_{p}^{s}$, $%
f_{a,d,l}(x)=p^{-3}f_{a,d}(l+px)$, $|f_{a,d,l}(x)|_{p}^{s}=|x|_{p}^{3s}$,$\
\nu (\overline{f_{a,d}})=p-2$, and $\delta (\overline{f_{a,d}})=1$, we
obtain \ that

\begin{eqnarray}
Z(s,f) &=&p^{-1}(p-3)+\frac{(1-p^{-1})p^{-1-s}}{1-p^{-1-s}}%
+p^{-1}(p-1)p^{-1-4s}  \notag \\
&&\,+p^{-1}(p-1)p^{-1-3s}+p^{-1}(p-2)p^{-2-8s}+\frac{(1-p^{-1})p^{-3-9s}}{%
1-p^{-1-s}}  \notag \\
&&\,+\frac{(1-p^{-1})p^{-3-11s}}{1-p^{-1-3s}}+p^{-2-6s}\int_{\mathbb{Z}%
_{p}}|f_{c,h}(x)|_{p}^{s}\mid dx\mid \text{.}  \label{step11}
\end{eqnarray}
Finally, by applying SPF to the integral \ involving $f_{c,h}(x)$ in (\ref
{step11}), we obtain that 
\begin{eqnarray}
Z(s,f) &=&p^{-1}(p-3)+\frac{(1-p^{-1})p^{-1-s}}{1-p^{-1-s}}%
+p^{-1}(p-1)p^{-1-4s}  \notag \\
&&\,+p^{-1}(p-1)p^{-1-3s}+p^{-1}(p-2)p^{-2-8s}+\frac{(1-p^{-1})p^{-3-9s}}{%
1-p^{-1-s}}  \notag \\
&&\,+\frac{(1-p^{-1})p^{-3-11s}}{1-p^{-1-3s}}+p^{-1}(p-2)p^{-2-6s}+\frac{%
(1-p^{-1})p^{-3-7s}}{1-p^{-1-s}}  \notag \\
&&+\frac{(1-p^{-1})p^{-3-8s}}{1-p^{-1-2s}}\text{.}  \label{step12}
\end{eqnarray}

\begin{remark}
\label{remark}If $\alpha =\frac{a}{b}\in \mathbb{Q}$, and $v_{p}(\alpha )<0$%
, then 
\begin{equation}
\mid x-\alpha \mid _{p}=\mid \alpha \mid _{p},\text{ for every }x\in \mathbb{%
Z}_{p}.  \label{const}
\end{equation}
On the other hand, a polynomial of the form 
\begin{equation*}
f(x)=\alpha _{0}\prod\limits_{i=1}^{r}(x-\alpha _{i})^{e_{i}}\in \mathbb{Q}%
[x],
\end{equation*}
can be decomposed as $f(x)=\alpha _{0}f_{-}(x)f_{+}(x)$, where 
\begin{equation}
f_{-}(x)=\prod\limits_{\{\alpha _{i}\mid v_{p}(\alpha _{i})<0\}}(x-\alpha
_{i})^{e_{i}},\text{ and \ }f_{+}(x)=\prod\limits_{\{\alpha _{i}\mid
v_{p}(\alpha _{i})\geqq 0\}}(x-\alpha _{i})^{e_{i}}.  \label{pol1}
\end{equation}
From (\ref{const}) and (\ref{pol1}) follow that \ 
\begin{equation*}
Z(s,f)=\mid \alpha _{0}\prod\limits_{\{\alpha _{i}\mid v_{p}(\alpha
_{i})<0\}}\alpha _{i}{}^{e_{i}}\mid _{p}^{s}Z(s,f_{+}).
\end{equation*}
Thus, from a computational point of view, we may assume without loss of
generality that all roots of $f(x)$ are $p-$adic integers.
\end{remark}

\section{Trees and $p$-adic Numbers}

The tree \ $U=U(p)$ of residue classes modulo powers of a given prime number 
$p$ is defined as \ follows. Consider the diagram 
\begin{equation*}
\{0\}=\mathbb{Z}/p^{0}\mathbb{Z}\underleftarrow{\text{ \ }\phi _{1}}\text{ \ 
}\mathbb{Z}/p^{1}\mathbb{Z}\text{ \ }\underleftarrow{\phi _{2}}\text{ \ }%
\mathbb{Z}/p^{2}\mathbb{Z}\text{ }\underleftarrow{\phi _{3}}\cdots
\end{equation*}
where $\phi _{l}$\ the \ are the natural homomorphisms. The vertices of \ $U$
are the elements of $\mathbb{Z}/p^{l}\mathbb{Z}$, for \ $l=0,1,2,\mathbb{%
\cdots }$, and the directed edges are $u\rightarrow v$ where $u\in \mathbb{Z}%
/p^{l}\mathbb{Z}$ and \ $\phi _{l}(u)=v$, for some $l>0$. Thus $U$ is a
rooted tree with root $\{0\}$. Exactly one directed edge emanates \ from
each vertex of $U$; except from the vertex $\{0\}$, from which no edge
emanates. In addition, every vertex is the end point of exactly $p$ directed
edges.

Given two \ vertices $u$, $v$ the notation $u>v$ will mean that there is a
sequence of vertices and edges of the form 
\begin{equation*}
u\rightarrow u^{(1)}\rightarrow \cdots \rightarrow u^{(m)}=v\text{.}
\end{equation*}
The notation $u\geqq v$ \ will mean that $u=v$ or $u>v$. The\textit{\ level} 
$l(u)$ of a vertex $u$ is $m$ if $u\in \mathbb{Z}/p^{m}\mathbb{Z}$. The 
\textit{valence} $Val(u)$ of a vertex $u$ is defined as the number of
directed edges whose end point is $u$.

A subtree, or simply a tree, is defined as a nonempty subset $T$ of vertices
of $U$, such that when $u\in T$ and $u>v$, then \ $v\in T$. Thus \ $T$
together with the directed edges $u\rightarrow v$, where $u,v\in T$ , is
again a tree with root $\{0\}$.

A tree $T$ is named a \textit{weighted tree}, if there exists a weight
function $W:T\rightarrow \mathbb{N}$. The value $W(u)$ is called the weight
of vertex $u$.

If \ $x\in \mathbb{Z}_{p}$, and $x_{l}$ denotes its residue class modulo $%
p^{l}$, then every vertex of $U$ is of the type $x_{l}$ with $l\in \mathbb{N}
$.

A \textit{stalk} is defined as a tree $K$ having at most one vertex at each
level. Thus a stalk is either finite, of the type 
\begin{equation*}
\{0\}\longleftarrow u^{(1)}\longleftarrow \cdots \longleftarrow u^{(l)},
\end{equation*}
or infinite, of the type 
\begin{equation*}
\{0\}\longleftarrow u^{(1)}\longleftarrow \cdots \text{.}
\end{equation*}
Clearly a finite stalk may be written as \ 
\begin{equation*}
\{0\}\longleftarrow x_{1}\longleftarrow \cdots \longleftarrow x_{l},
\end{equation*}
with $x\in \mathbb{Z}$, and infinite stalks as 
\begin{equation*}
\{0\}\longleftarrow x_{1}\longleftarrow x_{2}\longleftarrow \cdots ,
\end{equation*}
with $x\in \mathbb{Z}_{p}$. Thus there is a $1-1$ correspondence between \
infinite stalks and $p-$adic integers.

\subsection{Tree Attached to a Polynomial}

Let 
\begin{equation}
f(x)=\alpha _{0}\prod\limits_{i=1}^{r}(x-\alpha _{i})^{e_{i}}\in \mathbb{Q}%
[x]  \label{pol}
\end{equation}
be a non-constant polynomial, in one variable,\ of degree $d_{f}$, such that 
$v_{p}(\alpha _{i})\geqslant 0$, $i=1,2,\mathbb{\cdots },r$. \ We \
associate to $f(x)$ and a \ prime number $p$ the integer 
\begin{equation*}
l_{f}=\left\{ 
\begin{array}{cc}
1+\max \{v_{p}(\alpha _{i}-\alpha _{j})\mid i\neq j\text{, }1\leq i\text{, }%
j\leq r\}\text{,} & \text{if \ }r\geqq 2\text{;} \\ 
1\text{, } & \text{if \ \ }r=1\text{.}
\end{array}
\right.
\end{equation*}
We set 
\begin{equation*}
\alpha _{i}=a_{0,i}+a_{1,i\text{ }}p+\mathbb{\cdots }+a_{j,i}\text{ }p^{j}+%
\mathbb{\cdots }+a_{l_{f},i}\text{ }p^{l_{f}}\text{ mod }p^{l_{f}+1}\text{, }
\end{equation*}
$a_{j,i}\in \{0,1,\mathbb{\cdots },p-1\}$, $j=0,1,\mathbb{\cdots },l_{f}$, $%
i=1,2,\mathbb{\cdots },r$, for the $p-$adic expansion modulo $p^{l_{f}+1}$
of \ $\alpha _{i}$. We attach a weighted tree $T(f$, $l_{f})$ to $f$ as
follows: 
\begin{equation}
T(f,l_{f},p)=T(f,l_{f})=\bigcup\limits_{i=1}^{r}K(\alpha _{i},l_{f}),
\label{tree}
\end{equation}
where $K(\alpha _{i},l_{f})$\ denotes the stalk corresponding to the $p-$%
adic expansion of $\alpha _{i}$\ modulo $p^{l_{f}+1}$. Thus $T(f$, $l_{f})$\
is a rooted tree. We introduce a weight function on $T(f$, $l_{f}),$\ by
defining the weight of a vertex $u$ of level $m$ as 
\begin{equation}
W(u)=\left\{ 
\begin{array}{cc}
\sum\limits_{\{i\mid \alpha _{i}\equiv u\text{ mod }p^{m}\}}e_{i}, & \text{%
if }\ m\geqq 1\text{;} \\ 
0, & \text{if }\ m=0\text{.}
\end{array}
\right.  \label{weight}
\end{equation}

Given a vertex $u\in T(f,l_{f})$, we define the stalk generated by $u$\ to
be 
\begin{equation*}
B_{u}=\{v\in T(f,l_{f})\mid u\geqq v\}.
\end{equation*}
\ We associate a weight $W^{\ast }(B_{u})$ to \ $B_{u}$ as follows: 
\begin{equation}
W^{\ast }(B_{u})=\sum\limits_{v\in B_{u}}W(v).  \label{weightstalk}
\end{equation}

\subsection{Computation of Trees Attached to \ Polynomials}

Our next step is to show that a tree $T(f,l_{f})$ attached to a polynomial $%
f(x)$, of type (\ref{pol}), can be computed in polynomial time. There are
well known programming techniques \ to construct and manipulate trees and
forests (see e.g. \cite[Volume 1]{K}), for this reason, we shall focus on
showing \ that such computations can be carry out in polynomial time, and
set aside the\ implementation details of a particular algorithm for this
task. We shall include in the computation of \ $T(f,l_{f})$, the computation
of the weights of the stalks generated by its vertices; because all these
data will be used in the computation of the local zeta function of $f$.

\begin{proposition}
\label{prop4.1}The computation of a tree $T(f,l_{f})$ attached to a
polynomial $f(x)$, of type (\ref{pol}), \ from the \ $p-$adic expansions
modulo $p^{l_{f}+1}$\ of its roots 
\begin{equation*}
\alpha _{i}=a_{0,i}+a_{1,i\text{ }}p+\mathbb{\cdots }+a_{l_{f},i}p^{l_{f}}%
\text{ mod }p^{l_{f}+1}\text{ }
\end{equation*}
and multiplicities $e_{i}$, $i=1,2,\mathbb{\cdots },r$, involves \ $O(l_{f%
\text{ }}^{2}d_{f}^{3})$ arithmetic operations on integers with binary
length 
\begin{equation*}
O(\max \{\log p,\log (l_{f}d_{f})\}).
\end{equation*}
\end{proposition}

\begin{proof}
We assume that $T(f,l_{f})$ is finite set of the form 
\begin{equation}
T=\{\text{Level}_{0}\text{,}\cdots \text{,Level}_{j}\text{,}\cdots \text{%
,Level}_{l_{f}+1}\}\text{,}  \label{tipotree}
\end{equation}
where \ Level$_{j}$ represents the set of all vertices with level $j$. Each
Level$_{j}$ is a set of the form 
\begin{equation*}
Level_{j}=\{u_{j,1},\mathbb{\cdots },u_{j,i},\mathbb{\cdots },u_{j,m_{j}}\},
\end{equation*}
and each $u_{j,i}$ is a weighted vertex for every $i=1,\mathbb{\cdots },m_{j}
$. A \ weighted vertex $u_{j,i}$ is a set of the form 
\begin{equation*}
u_{j,i}=\{\text{ }W(u_{j,i})\text{, }Val(u_{j,i}),W^{\ast }(B_{u_{j,i}})\}%
\text{,}
\end{equation*}
where $W(u_{j,i})$ \ is the weight of $u_{j,i}$, $Val(u_{j,i})$ is its
valence, and $W^{\ast }(B_{ui})$ is \ the weight of \ stalk $B_{u_{j,i}}$.
The weight of the stalk generated by \ $u_{j,i}$ can be written as 
\begin{equation*}
W^{\ast }(B_{u_{j,i}})=\sum\limits_{v\in B_{u_{j,i}}}W(v).
\end{equation*}

For the computation of a vertex \ $u_{j,i}$\ of level $j$, we proceed as
follows. We put $I=\{1,2,\mathbb{\cdots },r\},$ and 
\begin{equation*}
M_{j}=\{\alpha _{i}\text{ mod }p^{j}\mid i\in I\}.
\end{equation*}
For each \ $0\leq j\leq l_{f}+1$, we compute a partition of $I$ of type 
\begin{equation}
I=\bigcup\limits_{i=1}^{l_{j}}I_{j,i},  \label{parti}
\end{equation}
such that 
\begin{equation*}
\alpha _{t}\text{ mod }p^{j}=\alpha _{s}\text{ mod }p^{j},
\end{equation*}
for every $t,s\in I_{j,i}$. Each subset $I_{j,i}$\ corresponds \ to a vertex 
$u_{j,i}$\ of level $j$. This \ computation requires\ $O(l_{f}r^{2})$
arithmetic operations on integers with binary length $O(\log p)$. Indeed,
the cost of computing a ``yes or no'' answer for the question: $\alpha _{t}$
mod $p^{j}=\alpha _{s}$ mod $p^{j}?$ is $O(j)$ comparisons of integers with
binary length $O(\log p)$. In the worst case, there are $r$ vectors $M_{j}$,
\ and the computation of partition (\ref{parti}), for a fixed $j$, involves
the comparison of $\alpha _{t}$ with $\alpha _{l}$ for $l=t+1,t+2,\mathbb{%
\cdots },r$. This computation requires\ $O(jr^{2})$ arithmetic operations on
integers with binary length $O(\log p)$. Since $j\leqq l_{f}+1$, the
computation of partition (\ref{parti}) requires\ $O(l_{f}r^{2})$ arithmetic
operations on integers with binary length $O(\log p)$.

The weight of the vertex \ $u_{j,i}$\ is given by the expression 
\begin{equation*}
W(u_{j,i})=\sum\limits_{k\in I_{j,i}}e_{k}.
\end{equation*}
Thus the computation of the weight of a vertex requires $O(r)$ additions of
integers with binary length $O(\log d_{f}r)$.

For the computation of the valence of \ $u_{j,i}$, we proceed as follows.
The valence of $u_{j,i}$ \ can be expressed as 
\begin{equation*}
Val(u_{j,i})=\text{Card}\{I_{j+1,l}\mid I_{j+1,l}\subseteq I_{j,i}\},
\end{equation*}
where $I_{j+1,l}$ runs through all possible sets that correspond to the \
vertices $u_{j+1,l}$, with level $j+1$. Thus the computation of $Val(u_{j,m})
$\ involves \ the computation of $\ $a ``yes or no'' answer for the question 
$I_{j+1,l}\subseteq I_{j,i}$? The computation of a\ ``yes or no'' answer \
involves $O(r)$\ comparisons of integers with binary length $O(\log r)$.
Therefore the computation of $Val(u_{j,i})$\ involves $O(r)$\ comparisons
and $O(r)$ additions of integers with binary length $O(\log r)$.

For the computation of the weight of \ $B_{u_{j,i}}$, we observe that $%
W^{\ast }(B_{u_{j,i}})$ is given by the formula 
\begin{equation*}
W^{\ast }(B_{u_{j,i}})=\sum\limits_{l=0}^{j-1}\sum\limits_{I_{j,i}\subseteq
I_{l,k}}W(I_{l,k}),
\end{equation*}
where $W(I_{l,k})=W(v_{l,k})$, and $v_{l,k}$\ is \ the vertex corresponding
to $I_{l,k}$. Thus the computation of $W^{\ast }(B_{u_{j,i}})$ involves $%
O(l_{f})$ additions of integers with binary length $O(\log (l_{f}$ $d_{f}))$%
, and $O(l_{f}$ $r)$ comparisons of integers with binary length $O(\log r)$.

From the above reasoning follows that the computation of a vertex of a tree $%
T(f,l_{f})$ \ involves at most $O(l_{f}$ $r^{2})$ arithmetic operations
(additions and comparisons) on integers with binary length $O(\max \{\log p,$
$\log ($ $l_{f}$ $d_{f})\})$. Finally, since the number of vertices of $%
T(f,l_{f})$ is at most $O(l_{f}$ $d_{f})$, it follows that \ the computation
of a tree of type $T(f,l_{f})$ involves $O(l_{f}^{2}$ $d_{f}^{3})$
arithmetic operations on integers with binary length $O(\max \{\log p,$ $%
\log ($ $l_{f}$ $d_{f})\})$.
\end{proof}

\section{Generating Functions and Trees}

In this section we attach to a weighted tree $T(f,l_{f})$ and a \ prime $p$
a generating function $G(s,T(f,l_{f}),p)\in \mathbb{Q}(p^{-s})$ defined \ as
follows.

We set 
\begin{equation*}
\mathcal{M}_{T(f,l_{f})}=\left\{ u\in T(f,l_{f})\left| 
\begin{array}{c}
W(u)=1\text{, and there no exists }v\in T(f,l_{f}) \\ 
\text{with }W(v)=1\text{, such that }u>v.
\end{array}
\right. \right\} ,
\end{equation*}

and

\begin{equation*}
L_{u}(p^{-s})=\left\{ 
\begin{array}{cc}
\frac{(1-p^{-1})p^{-l(u)-W^{\ast }(B_{u})s}}{(1-p^{-1-W(u)s})}\text{, \ if}
& l(u)=1+l_{f}\text{, and \ }W(u)\geqq 2\text{;} \\ 
&  \\ 
p^{-1}(p-Val(u))p^{-l(u)-W^{\ast }(B_{u})s}\text{, \ if} & 0\leqq l(u)\leqq
l_{f}\text{, and \ }W(u)\neq 1\text{;} \\ 
&  \\ 
\frac{(1-p^{-1})p^{-l(u)-W^{\ast }(B_{u})s}}{1-p^{-1-s}}\text{, \ \ if} & 
u\in \mathcal{M}_{T(f,l_{f})}\text{;} \\ 
&  \\ 
0\text{,\ \ \ if} & W(u)=1\text{, and }u\notin \mathcal{M}_{T(f,l_{f})}\text{%
.}
\end{array}
\right. 
\end{equation*}

With all the above notation, we define the generating function attached to $%
T(f,l_{f})$ and $p$\ \ as 
\begin{equation}
G(s,T(f,l_{f}),p)=\sum\limits_{u\in T(f,l_{f})}L_{u}(p^{-s})\text{.}
\label{Gfunction}
\end{equation}

Our next goal is to show that $G(s,T(f,l_{f}),p)=Z(s,f)$. The proof of this
fact requires the following preliminary result.

\begin{proposition}
The generating function attached to a tree $T(f,l_{f})$ and a prime $p$
satisfies 
\begin{eqnarray}
G(s,T(f,l_{f}),p) &=&p^{-1}\nu (\overline{f})+\delta (\overline{f})\frac{%
(1-p^{-1})p^{-1-s}}{(1-p^{-1-s})}  \notag \\
&&+\sum\limits_{\xi \in S}p^{-1-e_{\xi }s}G(s,T(f_{\xi },l_{f}-1),p).
\label{ident}
\end{eqnarray}
\end{proposition}

\begin{proof}
Let $A_{f}=\{u\in T(f,l_{f})\mid l(u)=1$, $W(u)=1\},$ and $B_{f}=\{u\in
T(f,l_{f})\mid l(u)=1$, $W(u)\geqq 2\}$. We have the following \ \ partition
for $T(f,l_{f}):$%
\begin{equation}
T(f,l_{f})=\{0\}\bigcup A_{f}\bigcup \left( \bigcup_{u\in B_{f}}T_{u}\right) 
\text{,}  \label{parti1}
\end{equation}
with 
\begin{equation*}
T_{u}=\{v\in T(f,l_{f})\mid v\geqq u\}.
\end{equation*}
Each $T_{u}$\ is a rooted tree with root $\{u\}$.\ \ From \ partition (\ref
{parti1}) and the definition of \ $G(s,T(f,l_{f}),p)$, it follows that 
\begin{eqnarray}
G(s,T(f,l_{f}),p) &=&p^{-1}\left( p-Val(\{0\})\right) +\text{Card}\{A_{f}\}%
\frac{(1-p^{-1})p^{-1-s}}{(1-p^{-1-s})}+  \notag \\
&&\sum\limits_{u\in B_{f}}G(s,T_{u}),  \label{ident0}
\end{eqnarray}
with $G(s,T_{u})=\sum\limits_{v\in T_{u}}L_{v}(p^{-s})$.

Since there exists a bijective correspondence between the roots of $%
\overline{f}(x)\equiv 0$ mod $p$ and the vertices $\ $of $T(f,l_{f})$\ with
level $1$, $\ \ \ \ \ $%
\begin{equation}
\ \ p-Val(\{0\})=\nu (\overline{f})\text{, and \ Card}\{A_{f}\}=\delta (%
\overline{f}).  \label{ident1}
\end{equation}

Now, if the vertex $u$ corresponds to the root \ $\overline{f}(\xi )\equiv 0$
mod $p$, $\ $\ then 
\begin{equation}
T_{u}=\left( \bigcup\limits_{\{\alpha _{i}\mid \alpha _{i}\equiv \xi \text{
mod }p\}}K(\alpha _{i},l_{f})\right) \setminus \{0\}.
\end{equation}
On the other hand, $\ $we have that 
\begin{equation}
T(f_{\xi },l_{f}-1)=\bigcup\limits_{\{\alpha _{i}\mid \alpha _{i}\equiv \xi 
\text{ mod }p\}}K(\frac{\alpha _{i-\xi }}{p},l_{f}-1).
\end{equation}
Now \ we remark that the map $\alpha _{i}\rightarrow \frac{\alpha _{i}-\xi }{%
p}$ induces a isomorphism between the trees $T_{u}$ and \ $T(f_{\xi
},l_{f}-1)$, that preserves the weights of the vertices; and \ thus we may
suppose that $T_{u}$ $=$\ $T(f_{\xi },l_{f}-1)$.\ The level function $l_{T}$
of $T(f_{\xi },l_{f}-1)$ is related to the level function $l_{T_{u}}$ of $%
T_{u}$ by means of the equality $l_{T}$ $-l_{T_{u}}=-1$. In addition, $B_{f}=
$ $S$, where $S$ is the subset of $\{0,1,\mathbb{\cdots },p-1\}\subseteq 
\mathbb{Z}_{p}$ whose reduction modulo $pZ_{p}$ is equal to the set of roots
of \ $\overline{f}(\xi )=0$ with multiplicity greater or equal than two. \
Therefore, it holds that \ 
\begin{equation}
G(s,T_{u})=p^{-1-e_{\xi }s}G(s,T(f_{\xi },l_{f}-1),p).  \label{ident2}
\end{equation}
The result follows from (\ref{ident0}) by the identities (\ref{ident1}) and (%
\ref{ident2}).
\end{proof}

\begin{lemma}
\label{lemma}Let $p$ be a fixed prime number and $v_{p}$ the corresponding $%
p-$adic valuation, and 
\begin{equation*}
f(x)=\alpha _{0}\prod\limits_{i=1}^{r}(x-\alpha _{i})^{e_{i}}\in \mathbb{Q}%
[x]\setminus \mathbb{Q},
\end{equation*}
a polynomial such that $v_{p}(\alpha _{i})\geqq 0$, for $i=1,\cdots ,r$.
Then 
\begin{equation*}
Z(s,f)=G(s,T(f,l_{f}),p).
\end{equation*}
\end{lemma}

\begin{proof}
We proceed by induction on $l_{f}$.

Case\ $l_{f}=1$

If $r=1$ the proof follows immediately, thus we may assume that $r\geqq 2$.
Since $l_{f}=1$, it holds that $v_{p}(\alpha _{i}-\alpha _{j})=0$, for every 
$i$, $j$, satisfying $i\neq j$, and thus \ $\overline{\alpha _{i}}\neq 
\overline{\alpha _{j}}$, if \ $i\neq j$. By applying SPF, we have that 
\begin{equation}
Z(s,f)=p^{-1}\nu (\overline{f})+\delta (\overline{f})\frac{(1-p^{-1})p^{-1-s}%
}{(1-p^{-1-s})}+\sum\limits_{\xi \in S}p^{-1-e_{\xi }s}\frac{(1-p^{-1})}{%
(1-p^{-1-e_{\xi }s})},
\end{equation}
where each $e_{\xi }=e_{j}\geqq 2$, for some $j$, and \ $\alpha _{j}=\xi
+p\beta _{j}$.

\ On the other hand, $T(f,l_{f})$ \ is a rooted tree with $r$ vertices $v_{j}
$, satisfying $l(v_{j})=1$, and $W(v_{j})=e_{j}$,\ for $j=1,\cdots ,r$.
These observations allow one to deduce that $Z(s,f)=G(s,T(f,l_{f}),p)$.

By induction hypothesis, we may assume that \ $Z(s,f)=G(s,T(f,l_{f}),p)$,
for every polynomial $f$ satisfying both the \ hypothesis of the lemma, and
the condition \ $1\leqq l_{f}\leqq k$, $k\in \mathbb{N}$.

Case\ $l_{f}=k+1$, $k\in \mathbb{N}$

Let $f(x)$ \ be \ a polynomial satisfying the lemma's \ hypothesis, and \ $%
l_{f}=k+1$, $k\geqq 1$. By applying SPF, we obtain that 
\begin{equation}
Z(s,f)=p^{-1}\nu (\overline{f})+\delta (\overline{f})\frac{(1-p^{-1})p^{-1-s}%
}{(1-p^{-1-s})}+\sum\limits_{\xi \in S}p^{-1-e_{\xi }s}\int \mid f_{\xi
}(x)\mid _{p}^{s}dx.  \label{eq1}
\end{equation}
Now, since $l_{f_{\xi }}=l_{f}-1$, for every \ $\xi \in S$, it follows from
the induction hypothesis applied to each \ $f_{\xi }(x)$ in (\ref{eq1}),
that 
\begin{equation}
Z(s,f)=p^{-1}\nu (\overline{f})+\delta (\overline{f})\frac{(1-p^{-1})p^{-1-s}%
}{(1-p^{-1-s})}+\sum\limits_{\xi \in S}p^{-1-e_{\xi }s}G(s,T(f_{\xi
},l_{f}-1),\text{ }p).  \label{eq3}
\end{equation}
Finally, from\ identity (\ref{ident}), and (\ref{eq3}), we conclude that 
\begin{equation}
Z(s,f)=G(s,T(f,l_{f}),p).  \label{eq4}
\end{equation}
\end{proof}

The following proposition gives a complexity estimate for the computation of 
$G(s,T(f,l_{f}),p)$.

\begin{proposition}
\label{prop4.2}The computation of the generating function \ 
\begin{equation*}
G(s,T(f,l_{f}),p)
\end{equation*}
from $T(f,l_{f})$, involves $O(l_{f}$ $d_{f})$ arithmetic operations on
integers with binary length \ $O(\max \{\log p,$ $\log ($ $l_{f}$ $d_{f})\})$%
.
\end{proposition}

\begin{proof}
This is a consequence of proposition \ref{prop4.1},\ and the definition of
generating function.
\end{proof}

\section{Computation of $p$-adic Expansions}

In this section we estimate the complexity of the steps $2$ and $3$ in \ the
algorithm Compute\_$Z(s,f)$.

\begin{proposition}
\label{prop4.3}Let 
\begin{equation*}
B=\max_{
\begin{array}{c}
{1\leq i,}\text{ }{j\leq r} \\ 
{i\neq j}
\end{array}
}{\{}\mid c_{j,i}\mid \text{,}\mid d_{j,i}\mid \text{ \ }{\mid }\text{ }%
\alpha _{j}-\alpha _{i}=\frac{c_{j,i}}{d_{j,i}}\text{, }c_{j,i}\text{, }%
d_{j,i}\in \mathbb{Z}\setminus \{0\}{\}}.
\end{equation*}
The computation of the integer $l_{f}$ involves $O(d_{f}^{2}\frac{\log B}{%
\log p})$ arithmetic operations on integers with binary length $O(\max
\{\log B,\log p\})$.
\end{proposition}

\begin{proof}
First, we observe that for $c\in \mathbb{Z}\setminus \{0\}$, \ the
computation of $v_{p}(c)$ \ involves $O(\frac{\log \mid c\mid }{\log p})$
divisions of integers of binary length $O(\max \{\log \mid c\mid ,$ $\log
p\})$. Thus the computation of $v_{p}(\frac{c}{d})=v_{p}(c)-v_{p}(d),$ \
involves $O(\frac{\max \{\log \mid c\mid ,\text{ }\log \mid d\mid \}}{\log p}%
)$ divisions and subtractions of integers with binary length 
\begin{equation*}
O(\max \{\log \mid c\mid ,\log \mid d\mid ,\log p\}).
\end{equation*}
From these observations follow that the computation of $v_{p}(\alpha
_{j}-\alpha _{i})$, $i\neq j$, $1\leq i$, $j\leq r$, involves $O(r^{2}\frac{%
\log B}{\log p})$ arithmetic operations on integers with binary length $%
O(\max \{\log B,\log p\})$. Finally, the computation of the maximum of the $%
v_{p}(\alpha _{j}-\alpha _{i})$, $i\neq j$, $1\leq i$, $j\leq r$, involves $%
O(\log r)$ comparisons of integers with binary length $O(\max \{\log B,\log
p\})$. Therefore the computation of the integer $l_{f}$ involves at most $%
O(d_{f}^{2}\frac{\log B}{\log p})$ arithmetic operations on integers with
binary length $O(\max \{\log B,\log p\})$.
\end{proof}

\begin{proposition}
Let \ $p$ be a fixed prime and $\gamma =\frac{c}{b}\in \mathbb{Q}$, with $%
c,b\in \mathbb{Z}\setminus \{0\}$, and $v_{p}(\gamma )\geqq 0$. The $p-$adic
expansion 
\begin{equation*}
\gamma =a_{0}+a_{1\text{ }}p+\mathbb{\cdots }+a_{j}p^{j}+\mathbb{\cdots }%
+a_{m}p^{m},
\end{equation*}
modulo $p^{m+1}$ involves $O(m+\log (\max \{\mid b\mid ,$ $p\}))$\
arithmetic operations on integers with binary length $O(\max \{\log \mid
c\mid ,\log \mid b\mid ,\log p\})$.
\end{proposition}

\begin{proof}
Let $y\in \{1,\cdots ,p-1\}$ be an integer \ such that $yb\equiv 1$ mod $p$.
This integer \ can be computed \ by means of the Euclidean algorithm in $%
O(\log (\max \{\mid b\mid ,p\}))$ \ arithmetic operations \ involving
integers of binary length $O(\max \{\log \mid b\mid ,\log p\})$ (cf. 
\cite[Volume 2, section 4.5.2]{K}).

We set $\gamma =\gamma _{0}=\frac{c}{b}$, $c_{0}=c$, and define \ $%
a_{0}\equiv yc$ mod $p$. With this notation, the $p-$adic digits $a_{i},i=1$,%
$\cdots $, $m$, can be computed recursively as follows: 
\begin{equation*}
\gamma _{i}=\frac{\frac{\left( c_{i-1}-a_{i-1}b\right) }{p}}{b}=\frac{c_{i}}{%
b},
\end{equation*}
\begin{equation*}
a_{i}=yc_{i}\text{ mod }p.
\end{equation*}
Thus the computation of the $p-$adic expansion of $\gamma $ \ needs \ $%
O(m+\log (\max \{\mid b\mid ,p\}))$\ arithmetic operations on integers with
binary length 
\begin{equation*}
O(\max \{\log \mid c\mid ,\text{ }\log \mid b\mid ,\text{ }\log p\}).
\end{equation*}
\end{proof}

\begin{corollary}
\label{cor4.4}Let $p$ be a fixed prime number and $v_{p}$ the corresponding $%
p-$adic valuation, and 
\begin{equation*}
f(x)=\alpha _{0}\prod\limits_{i=1}^{r}(x-\alpha _{i})^{e_{i}}\in \mathbb{Q}%
[x],
\end{equation*}
a non-constant polynomial such that $v_{p}(\alpha _{i})\geqq 0$, $i=1,\cdots
,r$. The computation of \ the $p-$adic expansions modulo $p^{l_{f}+1}$\ of
the roots $\alpha _{i}$, $i=1,2,\cdots ,r$, of $f(x)$\ \ involves $O(d_{f}$ $%
l_{f}+d_{f}$ $\log (\max \{B,p\}))$\ arithmetic operations \ on integers
with binary length \ $O(\max \{\log B,\log p\})$.
\end{corollary}

\begin{proof}
The corollary \ follows directly from the two previous propositions.
\end{proof}

\section{Computing local zeta functions of polynomials with splitting $%
\mathbb{Q}$}

In this section we prove the \ correctness of the algorithm Compute\_$Z(s,f)$
and estimate its complexity.

\begin{theorem}
\label{theoA}The algorithm Compute\_$Z(s,f)$ \ outputs\ the meromorphic
continuation of the Igusa local zeta function $Z(s,f)$ of a polynomial \ $%
f(x)\in \mathbb{Z}[x]$, in one variable, with splitting field \ $\mathbb{Q}$%
. The number of arithmetic operations needed by the algorithm is 
\begin{equation*}
O\left( d_{f}^{6}+d_{f}^{9}\log (\shortparallel f\shortparallel
)+l_{f}^{2}d_{f}^{3}+d_{f}^{2}\log \left( \max \{B,p\}\right) \right) ,
\end{equation*}
and the integers on which these operations \ are performed have a binary
length 
\begin{equation*}
O\left( \max \{\log p\text{, }\log l_{f}d_{f}\text{, }\log B\text{, }%
d_{f}^{3}+d_{f}^{2}\log (\shortparallel f\shortparallel )\}\right) \text{.}
\end{equation*}
\end{theorem}

\begin{proof}
By remark (\ref{remark}), we may assume without loss of generality that 
\begin{equation*}
f(x)=\alpha _{0}\prod\limits_{i=1}^{r}(x-\alpha _{i})^{e_{i}}\in \mathbb{Q}%
[x]\setminus \mathbb{Q},
\end{equation*}
with $v_{p}(\alpha _{i})\geqq 0$, $i=1,\cdots ,r$. The correctness of the
algorithm follows from lemma \ \ref{lemma}. The complexity estimates are
obtained as follows: the number of arithmetic operations \ needed in the
steps 2 (cf. proposition \ref{prop4.3}), 3 (cf. \ corollary \ref{cor4.4}), 4
(cf. proposition \ref{prop4.1}), 5 (proposition\ \ref{prop4.2}), and 6 is at
most 
\begin{equation*}
O\left( l_{f}^{2}d_{f}^{3}+d_{f}^{2}\log \left( \max \{B,p\}\right) \right) 
\text{;}
\end{equation*}
and these operations are performed on integers \ whose binary length is at
most 
\begin{equation*}
O\left( \max \{\log p\text{, }\log l_{f}d_{f}\text{,}\log B\}\right) .
\end{equation*}
The estimates for the whole algorithm \ follow from the above estimates and
those of \ the factoring algorithm by \ A. K. Lenstra, H. Lenstra and L. Lov%
\'{a}sz (see theorem 3.6 of \cite{L-L-L}).
\end{proof}

\section{Stream Ciphers and Poincar\'{e} series}

There is a natural connection between Poincar\'{e} series and stream
ciphers. In order to explain this relation, we recall some basic facts about
stream ciphers \ \cite{R}. \ Let $\mathbb{F}_{p^{n}}$ be a finite field with 
$p^{n}$ elements, with $p$ a prime number. \ For any integer $r>0$ and $r$\
fixed elements $q_{i}\in \mathbb{F}_{p^{n}}$, $i=1,\mathbb{\cdots },r$\
(called taps), a Linear Feedback Shift Register, abbreviated LFSR, of length 
$r$ consists of $r$ cells with initial contents $\left\{ a_{i}\in \mathbb{F}%
_{p^{n}}\mid i=1,\mathbb{\cdots },r\right\} $. For any $n\geqslant r$, if
the current state is $(a_{n-1},\mathbb{\cdots },a_{n-r})$, then $a_{n}$ is
determined by the linear recurrence relation 
\begin{equation*}
a_{n}=-\sum\limits_{i=1}^{r}a_{n-i}q_{i}.
\end{equation*}
The device outputs the rightmost element $a_{n-r}$, shifts all the cells one
unit right, and feeds $a_{n}$ back to the leftmost cell. 

Any configuration of the $r$ cells \ forms a state of the LSFR. If $%
q_{r}\neq 0$, the following polynomial \ $q(x)\in \mathbb{F}_{p^{n}}[x]$ of
degree $r$\ appears in the analysis of LFSRs: 
\begin{equation*}
q(x)=q_{0}+q_{1}x+\mathbb{\cdots }.+q_{r}x^{r}\text{ \ with }q_{0}=-1\text{.}
\end{equation*}
This polynomial is called the connection polynomial. An infinite sequence $%
A=\left\{ a_{i}\in \mathbb{F}_{p^{n}}\mid i\in \mathbb{N}\right\} $ has
period $T$ if for any $i\geqslant 0$, \ $a_{i+T}=a_{i}$. Such a sequence is
called periodic. If this is only true for $i$ greater than some index $i_{0}$%
, then the sequence is called eventually periodic. The following facts about
an LFSR of length $r$ are well-known \cite{R}.

\begin{enumerate}
\item  There are only finitely many possible states, and the state with all
the cells zero will produce a $0-$sequence. The output sequence is
eventually periodic and the maximal period is $p^{nr}-1$.

\item  The Poincar\'{e} series $\ g(x)=\sum\limits_{i=0}^{\infty }a_{i}x^{i}$
associated with the output sequence is called the generating function of the
sequence. It is a rational function over $\mathbb{F}_{p^{n}}$\ of the form $%
g(x)=\frac{L(x)}{R(x)}$, with $L(x)$, $R(X)$ $\in \mathbb{F}_{p^{n}}[x]$, $%
deg(R(X))<r$. The output sequence is strictly periodic if and only if $%
deg(L(X))<deg(R(X))$.

\item  There is a one-to-one correspondence between LFSRs of length $r$ with 
$q_{r}\neq 0$ and rational functions $\frac{L(x)}{R(x)}$ with $deg(R(X))=r$
and $deg(L(X))<r$.
\end{enumerate}

We set $\mathbb{F}_{p^{n}}(x)$ for the field of rational functions over $%
\mathbb{F}_{p^{n}}$, and $N^{\infty }(\mathbb{F}_{p^{n}})$ for the set of
sequences \ of the form $\{b_{0},\mathbb{\cdots },b_{u}\}$, $b_{i}\in 
\mathbb{F}_{p^{n}}$, $\ 0\leq i\leq u$, $u\in \mathbb{N}$. From the above
considerations, it is possible to identify an LFSR with a function $F_{u}$, $%
u\in \mathbb{N}$, defined as follows: 
\begin{equation}
\begin{array}{cccc}
F_{u}: & \mathbb{F}_{p^{n}}(x) & \rightarrow & N^{\infty }(\mathbb{F}%
_{p^{n}}) \\ 
& \sum\limits_{i=0}^{\infty }a_{i}x^{i} & \rightarrow & \{a_{0},\mathbb{%
\cdots },a_{u}\}.
\end{array}
\label{A}
\end{equation}
We set 
\begin{equation*}
\mathcal{H}=\{H(t,f)\mid f(x)\in \mathbb{Z}[x]\text{, in \ one variable,
with splitting field }\mathbb{Q}\},
\end{equation*}
and $N^{\infty }$ $\left( \mathbb{Z}\right) $\ for the set \ of finite
sequences of integers. Also, for each $u\in \mathbb{N}$, and a prime number $%
p$, we define 
\begin{equation}
\begin{array}{cccc}
F_{u,p}: & \mathcal{H} & {\ \rightarrow } & \mathbb{N}^{\infty }\left( 
\mathbb{Z}\right) \\ 
& H(t,f) & {\ \rightarrow } & \{N_{0}(f,p),N_{1}(f,p),\mathbb{\cdots }%
,N_{u}(f,p)\}\text{.}
\end{array}
\label{B}
\end{equation}
Thus the mappings $F_{u,p}$ can be seen as LFSRs, or stream ciphers, over $%
\mathbb{Z}$. If we replace each $N_{u}(f,p)$ by its binary representation,
then the $\ F_{u,p}$ are LFSRs. For practical purposes it is necessary that $%
F_{u,p}$ \ can be computed efficiently, i.e., in polynomial time. With the
above notation our second result is the following.

\begin{theorem}
\label{theoB} For every $H(t,f)$\ $\in $ $\mathcal{H}$, the computation of \ 
$F_{u,p}(H(t,f))$ involves $O(u^{2}d_{f}l_{f})$ arithmetic operations, and
the integers on which these operations are performed have binary length 
\begin{equation*}
O(\max \{\text{ }(l_{f}+u)\log p,\text{ }\log (d_{f}l_{f})\}).
\end{equation*}
\end{theorem}

The proof of this theorem will be given at the end of this section. This
proof requires some preliminary results. We set $t=q^{-s}$, and 
\begin{equation*}
Z(s,f)=Z(t,f)=\sum\limits_{m=0}^{\infty }c_{m}(f,p)t^{m},
\end{equation*}
with $c_{m}(f,p)=vol(\{x\in \mathbb{Z}_{p}\mid v_{p}(f(x))=m\})$.

\begin{proposition}
Let $f(x)\in \mathbb{Z}[x]\setminus \mathbb{Z}$ be a polynomial in one
variable and $p$ a prime number. The following formula holds for $N_{n}(f,p)$%
: 
\begin{equation}
N_{n}(f,p)=\left\{ 
\begin{array}{ccc}
1\text{,} & \text{if} & n=0\text{;} \\ 
p^{n}\left( 1-\sum\limits_{j=1}^{n}c_{j-1}(f,p)\right) \text{,} & \text{if}
& n\geqslant 1\text{.}
\end{array}
\right.  \label{For1}
\end{equation}
\end{proposition}

\begin{proof}
The result follows by comparing \ the coefficient of $t^{n}$ of the series 
\begin{equation*}
\sum\limits_{n=0}^{\infty }\frac{N_{n}(f,p)}{p^{n}}t^{n}\text{ \ \ and \ }%
\sum\limits_{n=0}^{\infty }d_{n}t^{n}\text{,}
\end{equation*}
in the following equality : 
\begin{equation*}
H(t,f)=\sum\limits_{n=0}^{\infty }\frac{N_{n}(f,p)}{p^{n}}t^{n}=\frac{%
1-t\left( \sum\limits_{m=0}^{\infty }c_{m}(f,p)t^{m}\right) }{1-t}%
=\sum\limits_{n=0}^{\infty }d_{n}t^{n}.
\end{equation*}
\end{proof}

We associate to each $u\in T(f,l_{f})$, and $j\in \mathbb{N}$, a rational
integer $a_{j}(u)$ defined as follows:

\begin{equation}
a_{j}(u)=\left\{ 
\begin{array}{cc}
\frac{(p-1)}{p^{l(u)+1+y(u)}}\text{, \ if} & 
\begin{array}{c}
l(u)=1+l_{f}\text{, }W(u)\geqq 2\text{, }j=W^{\ast }(B_{u})+y(u)W(u)\text{,}
\\ 
\text{for some }y(u)\in \mathbb{N}\text{;}
\end{array}
\\ 
&  \\ 
\frac{(p-Val(u))}{p^{l(u)+1}}\text{, \ if} & 0\leqq l(u)\leqq l_{f}\text{, }%
W(u)\neq 1\text{, }j=W^{\ast }(B_{u})\text{;} \\ 
&  \\ 
\frac{(p-1)}{p^{l(u)+1+y(u)}}\text{, \ if} & 
\begin{array}{c}
u\in \mathcal{M}_{T(f,l_{f})},j=W^{\ast }(B_{u})+y(u)\text{,} \\ 
\text{for some }y(u)\in \mathbb{N}\text{;}
\end{array}
\\ 
&  \\ 
0\text{, \ \ \ \ \ \ \ \ \ \ \ if} & W(u)=1\text{, and }u\notin \mathcal{M}%
_{T(f,l_{f})}\text{;} \\ 
&  \\ 
0\text{, \ \ \ \ \ \ \ \ \ \ \ \ } & \text{in other cases.}
\end{array}
\right.  \label{For2}
\end{equation}

\begin{proposition}
Let $f(x)\in \mathbb{Z}[x]\setminus \mathbb{Z}$ be a polynomial in one
variable, with splitting field $\mathbb{Q}$, and $p$ a prime number. The
following formula holds: 
\begin{equation}
c_{j}(f,p)=\sum\limits_{u\in T(f,l_{f})}a_{j}(u),\text{ }j\geqq 0.
\label{For4}
\end{equation}
\end{proposition}

\begin{proof}
As a consequence of lemma (\ref{lemma}), we have the following identity: 
\begin{equation}
Z(t,f)=\sum\limits_{u\in T(f,l_{f})}L_{u}(t),  \label{For4a}
\end{equation}
with 
\begin{equation}
L_{u}(t)=\left\{ 
\begin{array}{cc}
\frac{(p-1)t^{W^{\ast }(B_{u})}}{p^{l(u)+1}(1-p^{-1}t^{\text{ }W(u)})}\text{%
, \ if} & l(u)=1+l_{f}\text{, }W(u)\geqq 2\text{;} \\ 
&  \\ 
\frac{(p-Val(u))}{p^{l(u)+1}}t^{W^{\ast }(B_{u})}\text{, \ if} & 0\leqq
l(u)\leqq l_{f}\text{, }W(u)\neq 1\text{;} \\ 
&  \\ 
\frac{(p-1)t^{W^{\ast }(B_{u})}}{p^{l(u)+1}(1-p^{-1}t)}\text{, \ if} & u\in 
\mathcal{M}_{T(f,l_{f})}\text{;} \\ 
&  \\ 
0\text{, \ \ \ \ \ if} & W(u)=1\text{, and }u\notin \mathcal{M}_{T(f,l_{f})}.
\end{array}
\right.  \label{For5}
\end{equation}

The result follows by comparing the coefficient of $t^{j}$ in the series $%
Z(t,f)=\sum\limits_{m=0}^{\infty }c_{m}(f,p)t^{m}$, and $Z(t,f)=\sum%
\limits_{u\in T(f,l_{f})}L_{u}(t)$.
\end{proof}

\begin{proposition}
\label{prop8}Let $f(x)\in \mathbb{Z}[x]\setminus \mathbb{Z}$ be a polynomial
in one variable, with splitting field $\mathbb{Q}$, and $p$ a prime number.

\begin{enumerate}
\item  The computation \ of \ $N_{n}(f,p)$, $n\geqq 1$, \ from \ the $%
c_{j-1}(f,p)$, $j=1,\mathbb{\cdots },n$, involves \ $O(n)$ \ \ arithmetic
operations on integers with binary length \ $O(n\log p)$.

\item  The computation \ of \ $c_{j}(f,p)$, $j\geqq 0$, \ from \ $Z(t,f)$,
involves \ $O(d_{f}l_{f})$ \ \ arithmetic operations on integers \ with
binary length \ 
\begin{equation*}
\mathit{O}(\max \{(j+l_{f})\log p,\log p,\log (d_{f}l_{f})\}).
\end{equation*}

\item  The computation \ of any $N_{n}(f,p)$, $n\geqq 1$, from $Z(t,f)$,
involves $O(nd_{f}l_{f})$ \ arithmetic operations on integers with binary
length \ 
\begin{equation*}
\mathit{O}(\max \{(n+l_{f})\log p,\log (d_{f}l_{f})\}).
\end{equation*}
\end{enumerate}
\end{proposition}

\begin{proof}
(1) By (\ref{For2})\ and (\ref{For4}), $c_{j}(f,p)=\frac{v_{j}}{p^{m_{j}}}$, 
$v_{j},m_{j}\in \mathbb{N}$. In addition, 
\begin{equation*}
c_{j-1}(f,p)=p^{-j+1}N_{j-1}(f,p)-p^{-j}N_{j}(f,p).
\end{equation*}
Thus $p^{n}c_{j-1}(f,p)\in \mathbb{N}$, for $\ j=1,\mathbb{\cdots },n$, and
\ $m_{j}\leq $\ $n,$ for $\ j=1,\mathbb{\cdots },n$. From (\ref{For1}), it
follows that $\ $%
\begin{equation}
N_{n}(f,p)=p^{n}-\sum\limits_{j=1}^{n}p^{n}c_{j-1}(f,p),\text{ }n\geqslant 1.
\label{For3}
\end{equation}
The above formula implies that the computation of \ $N_{n}(f,p)$, $n\geqq 1$%
, \ from \ the $c_{j-1}(f,p)$, $j=1,\mathbb{\cdots },n$, involves \ $O(n)$ \
\ arithmetic operations on integers with binary length \ $O(n\log p)$.

(2) The \ computation of $a_{j}(u)$ from $L_{u}(t)$ (i.e. from $Z(t,f)$, cf.
(\ref{For4a}))\ involves $O(1)$ arithmetic operations (cf. (\ref{For2}), (%
\ref{For5})) \ on integers of binary length $O(\max \{\log p$, $\log
(d_{f}l_{f})\})$. Indeed, since the numbers $l(u)$, $W^{\ast }(B_{u})$, $%
W(u) $, $u\in T(f,l_{f})$ are involved in this computation, we know \ by
proposition \ref{prop4.1} that their binary length is bounded by $O(\max
\{\log p$, $\log (d_{f}l_{f})\})$.

The cost of computing \ $c_{j}(f,p)$ from $L_{u}(t)$, $u\in T(f,l_{f})$
(i.e. from $Z(t,f)$) is bounded by the number of vertices of $T(f,l_{f})$ \
multiplied by an upper bound \ for the cost of computing $a_{j}(u)$ from $%
L_{u}(t)$, for any $j$, and $u$ (cf. (\ref{For4})). Therefore, \ from the
previous discussion the cost of computing $c_{j}(f,p)$ from $Z(t,f)$ is
bounded by $O(d_{f}$ $l_{f})$ \ \ arithmetic operations. These arithmetic
operations are performed on integers of binary length bounded \ by $O(\max
\{(j+l_{f})\log p,\log p,$ $\log (d_{f}l_{f})\}).$ Indeed, the binary
lengths of the numerator and the denominator of \ $a_{j}(u)+a_{j}(u^{^{%
\prime }})$, $u$, $u^{^{\prime }}\in T(f,l_{f})$ are bounded by $%
(l_{f}+1+j)\log p$ (cf. (\ref{For2})). Thus, the mentioned arithmetic
operations \ for calculating $c_{j}(f,p)$ from\ $L_{u}(t)$ are performed on
integers whose binary length is bounded by \ $O(\max \{(j+l_{f})\log p,\log
p,$ $\log (d_{f}l_{f})\})$.

(3) The third part \ follows the first and second parts by (\ref{For3}).
\end{proof}

\subsection{Proof of Theorem \ref{theoB}}

The theorem follows from proposition \ref{prop8} (3).

\end{document}